\documentclass[review]{elsarticle}

\usepackage{lineno,hyperref}
\usepackage{amsmath}
\usepackage{amssymb}
\usepackage{mathtools}
\usepackage{mathabx}
\usepackage{hyperref}
\newtheorem{theorem}{Theorem}
\newtheorem{lemma}{Lemma}
\newtheorem{corollary}{Corollary}
\newtheorem{conjecture}{Conjecture}
\newproof{pf}{Proof}
\newdefinition{rmk}{Remark}
\newdefinition{definition}{Definition}
\usepackage{enumitem}
\usepackage{tikz}
\modulolinenumbers[5]

\journal{---}

\begin{document}

\begin{frontmatter}



\title{Universal Bounds for Spreading on Networks}


\author[mymainaddress]{Gadi Fibich\corref{mycorrespondingauthor}}
\cortext[mycorrespondingauthor]{Corresponding author}
\ead{fibich@tau.ac.il}

\author[mymainaddress]{Tomer Levin}

\address[mymainaddress]{Department of Applied Mathematics, School of Mathematical Sciences, Tel Aviv University, Tel Aviv 6997801, Israel}

\begin{abstract}
	Spreading (diffusion) of innovations is a stochastic process on social networks. When the key driving mechanism is peer effects (word of mouth), the rate at which the aggregate adoption level increases with time depends strongly on the network structure. In many applications, however, the  network structure is unknown. To estimate the aggregate adoption level for such innovations, we show that the two networks that correspond to the slowest and fastest adoption levels are a homogeneous two-node network and a homogeneous infinite  complete network, respectively.  Solving the stochastic Bass model on these two networks yields explicit lower and upper bounds for the adoption level on any network. These bounds are tight, and they also hold for the individual adoption probabilities of nodes. The gap between the lower and upper bounds increases monotonically with the ratio of the rates of internal and external influences.
\end{abstract}


\end{frontmatter}

\section{Intorduction}
Diffusion (spreading) in networks is an active research area  
in mathematics, economics, management science, physics, biology, computer science, and social sciences,
and it concerns the spreading of diseases, computer viruses, rumors, information, opinions, technologies,
innovations, and more~\citep{Albert-00,Anderson-92,Jackson-08,Pastor-Satorras-01,Strang-98}. 
In marketing, diffusion of new products is a classical problem~\citep{Mahajan-93}.

The first mathematical model of diffusion of new products
was introduced by Bass~\cite{Bass-69}.  In this model, individuals adopt a new product because
of {\em external influences} by mass media and commercials, and because of {\em internal influences} (peer effect, word-of-mouth)
by individuals who have already adopted the product.
Let~$f$ denote the {\em adoption level} (fraction of adopters) in the population at time~$t$. Then according to the Bass model, 
	\begin{equation}
		\label{eq:homogeneous_Bass_frac}
		f'(t)= \left(1-f\right)\left(p+qf\right), \quad t>0, \qquad f(0)=0.
	\end{equation}
	Thus, the $1-f$ potential adopters adopt due to external influences at the constant rate of~$p$, and due to internal influences at the rate of~$qf$, which is proportional to the fraction of adopters.
Equation~\eqref{eq:homogeneous_Bass_frac} can be solved explicitly, yielding the S-shaped Bass formula~\citep{Bass-69}
	\begin{equation}
		\label{eq:Bass_sol}
		f_{\rm Bass}(t)=\frac{1-e^{-(p+q)t}}{1+\frac{q}{p}e^{-(p+q)t}}.
	\end{equation}

The Bass model~\eqref{eq:homogeneous_Bass_frac} inspired a
huge body of theoretical and empirical research; in 2004 it  
was selected as one of the 10 most-cited papers in the 50-year
history of Management Science~\citep{Hopp-04}. Initially, this research was carried out using {\em compartmental Bass models}, such as~\eqref{eq:homogeneous_Bass_frac}, in which the population is divided into several  compartments (e.g., nonadopters and adopters), and the transition
rates between compartments are given by deterministic ordinary differential equations.  {\it Compartmental Bass models, therefore, 
	implicitly assume  that the underlying social network is a homogeneous complete graph}, i.e., that
all individuals within the population are equally likely to influence each other.

In order not to make these assumptions, in more recent studies
diffusion of new products has been studied 
using  Bass models on networks, for the stochastic adoption decision of each individual~\citep{Delre-07, OR-10, Garber-04, Garcia-05, Muller-Peres-19, Rand-11}. 
These {\em agent-based models} allow for implementing a network structure, so that individuals are only influenced by adopters who are also their peers. For example, it has been suggested that social networks have a small-worlds~\citep{watts1998collective} or a scale-free structure~\citep{barabasi1999emergence}. 
In large-scale online social networks,   40\% of the links
were found to be within a 100 km radius~\citep{scellato2011socio}.
In the diffusion of solar panels, the key predictor of a new solar installation is having a neighbor who already installed one~\citep{Bollinger-12,Graziano-15}, and so  the relevant network is a two-dimensional Cartesian grid.  Bass models on networks also allow for heterogeneity among individuals~\cite{PNAS-12,fibich2022diffusion, GLM-01}.

   The  effects of 
	{\em various network  characteristics}  (average degree, clustering, \dots)  
on the diffusion were studied numerically using agent-based simulations, see e.g.,~\cite{Goldenberg-10,Goldenberg-08,Peres-14}. For example, it was found that growth is particularly effective in networks that demonstrate cohesion (strong mutual influence among its members), connectedness (high number of ties),
and conciseness (low redundancy) \cite{Muller-Peres-19}.

{\it Explicit expressions} for the expected adoption level~$f(t)$ in the  Bass model were only obtained for a few networks.
Niu \cite{Niu-02} explicitly computed the expected adoption level~$f_{\rm complete}(t;M)$ on complete homogeneous networks with $M$~nodes, and showed that
$\lim_{M \to \infty} f_{\rm complete}(t;M) = f_{\rm Bass}(t)$, see  Theorem~\ref{thrm:Niu} below. 
Fibich and Gibori \cite{OR-10} explicitly computed the expected adoption level~$f_{\rm circle}(t;M)$ on  homogeneous circles with $M$~nodes. They showed that the adoption level on the infinite circle,
 denoted by $f_{\rm 1D}(t):=\lim_{M \to \infty} f_{\rm circle}(t;M)$,
  is given by $f_{\rm 1D}(t)=1-e^{-(p+q)t+q\frac{1-e^{-pt}}{p}}$.

For most networks,  explicit expressions for~$f(t)$ are not available. Moreover,
in many applications, the network structure or even its characteristics are not known. Hence, {\em it is important, for both theoretical and practical considerations, to obtain explicit lower and upper bounds for the expected adoption level~$f(t)$.} 

In \cite{OR-10} it was conjectured that since circular and complete networks are the ``least-connected" and the ``most-connected" networks, the adoption level on any infinite network should be bounded from below by that on the infinite circle, and from above by that on the infinite complete network, i.e., that
$
f_{\rm 1D}(t)\leq f(t) \leq f_{\rm Bass}(t). 
$
So far, this conjecture has remained open.

In this study, we settle this conjecture. We 
prove that $f(t)\leq f_{\rm Bass}(t)$ for any finite or infinite network.
Thus, as was conjectured in~\cite{OR-10},  $f_{\rm Bass}$ is a universal upper bound for the adoption level.
Moreover, this upper bound is tight, and is strict for non-complete networks. The tight universal upper bound for 
the individual adoption probabilities of nodes
(i.e., for the probability of any node to adopt the product before time~$t$)
is also given by~$f_{\rm Bass}$.
 
The universal lower bound for~$f(t)$ on general finite or infinite networks, however, is not~$f_{\rm 1D}$. Rather, we prove that $f(t)\geq f_{M=2}^{\rm hom}$ 
for any network,
where~$f^{\rm hom}_{M=2}:= 1-e^{-pt}\frac{qe^{-pt}-pe^{-qt}}{q-p}$ is the expected adoption level on a homogeneous two-node network. 
This universal lower bound is also tight, and it also holds for the individual adoption probabilities of nodes. Thus, the conjecture from~\cite{OR-10} that~$f_{\rm 1D}$ is a universal lower bound for all infinite networks is wrong
(note, however, that for any~$D \ge 1$, $f_{\rm 1D}$ is the tight lower bound for the adoption level~$f_{\rm D}(t)$ on infinite D-dimensional Cartesian network 
where each node is connected to its 2D nearest neighbors with edges of weight $\frac{q}{2D}$, 
see~\cite{Funnel2022} for more details).

Let us motivate the ``success" of the conjecture from~\cite{OR-10} regarding the upper bound, and its ``failure" regarding the lower bound. As noted, the compartmental Bass model~\eqref{eq:homogeneous_Bass_frac} corresponds to a complete network, which is indeed the ``most-connected'' network, in the sense that each node can be  directly influenced by all other nodes. A one-sided circle, where each node can only influenced by the node to its left, however, is not the ``least-connected" network. This is because each node is also indirectly influenced by all other nodes. Rather, the ``least-connected" network is a  collection of disjoint pairs of nodes, where each node can be directly influenced by the other node in the pair, but cannot be indirectly influenced by any other node.

To quantify the influence of the social-network structure on the adoption level of new products, we study the {\em size of the gap between the lower and upper bounds}.
The gap size is a monotonically-increasing function of the 
ratio~$\frac{q}{p}$ of the rates of internal and external influences. 
For products that spread predominantly through word of mouth,
we obtain an explicit approximation for the gap size. This explicit approximation shows that the network structure indeed has a large influence on the adoption level of such products.

The paper is organized as follows. Section~\ref{sec:discrete_Bass} presents the   Bass model on a general network. Section~\ref{sec:main_results} presents the main results of this paper on the universal  lower and upper bounds.
Section~\ref{sec:gap} considers the size of the gap between the lower and upper bounds.
Section~\ref{sec:open} lists some open research problems. 
 The detailed proofs are given in Section~\ref{sec:proofs}.

\section{Bass model on networks}
	\label{sec:discrete_Bass}
	
	We begin by introducing the  Bass model on a general heterogeneous network. A new product is introduced at time $t=0$ to a network with 
	$M$~individuals, denoted by~${\cal M}:=\{1,\ldots,M\}$, where $M$ can be finite or infinite. We denote by~$X_j(t)$ the {\em state} of individual $j$ at time $t$, so that 
	\begin{equation*}
		X_j(t)=\begin{cases}
			1, \quad {\rm if}\ j\ {\rm is\ an\ adopter\ of\ the\ product\ at\ time}\ t,\\
			0, \quad {\rm otherwise.}
		\end{cases}
	\end{equation*} 
	Since the product is new, all individuals are initially nonadopters, i.e., 
	\begin{subequations}
	\label{eq:dbm}
	\begin{equation}
		\label{eq:general_initial}
		X_j(0)=0, \qquad j\in {\cal M}.
	\end{equation}
 The underlying social network is represented by a {\em weighted directed} graph, such that if there is an  edge from $k$ to $j$, the rate of {\em internal influence} of adopter~$k$ on nonadopter~$j$ to adopt is~$q_{k,j}>0$, and~$q_{k,j}=0$ if there is no edge from~$k$ to~$j$.
 	The edges and influence rates are {\em not} assumed to be symmetric, i.e., 
 $q_{k,j}$~may be different from~$q_{j,k}$. 
Since nonadopters do not self-influence to adopt,
$$
q_{j,j}\equiv 0, \qquad j\in {\cal M}.
$$	
In contrast to similar models in epidemiology on networks~\citep{kiss2017mathematics}, such as the Susceptible Infected (SI) model, $j$ also experiences {\em external influences}  to adopt by mass media and commercials, at a constant rate of~$p_j>0$. 
Internal and external influences are assumed to be additive. Thus,  
{\em the adoption time~$T_j$ of nonadopter~$j$ is exponentially distributed} at the rate of
\begin{equation}
\label{eq:lambda_j}
\lambda_j(t):=p_j+\sum_{k\in {\cal M}}q_{k,j}X_k(t), \qquad j\in {\cal M}, \quad t>0,
\end{equation}
which increases whenever $k$~adopts and $q_{k,j}>0$. Finally, it is assumed that once an individual adopts the product, she or he remains an adopter for all time.
Therefore, the stochastic adoption of~$j\in {\cal M}$ in the time interval~$(t,t+\Delta t)$ as~$\Delta t\rightarrow 0$ is given by
	\begin{equation}
		\label{eq:general_model}
		\mathbb{P}(X_j(t+\Delta t)=1\,| \, {\bf X}(t))=
		\begin{cases}
			\mbox{} \qquad \qquad  1,\hfill &\quad {\rm if}\ X_j(t)=1,\\
			\left(p_j+\sum\limits_{\substack{k \in {\cal M}}}q_{k,j}X_{k}(t)\right)\Delta t, &\quad {\rm if}\ X_j(t)=0,
		\end{cases}
	\end{equation}
\end{subequations}
	where ${\bf X}(t) := \left(X_1(t),\ldots, X_M(t)\right)$ is the {\em state of the network} at time~$t$. Note that the time variable is continuous.

	The maximal rate of internal influences that can be exerted on node~$j$ (which is when all its neighbors/peers  are adopters) is 
	\begin{subequations}
\begin{equation}
		\label{eq:q_on_node}
		q_j:=\sum_{\substack{k \in {\cal M}}}q_{k,j}.
	\end{equation}
	For simplicity, we assume that each node can be influenced by at least one node, i.e.,
	\begin{equation}
	\label{eq:q_j>0}
	q_j>0, \qquad j\in {\cal M}.
	\end{equation}
\end{subequations}	
 We do not assume, however, that the network  only consists of  a single connected component. 
 	The underlying network of the  Bass model~\eqref{eq:dbm} is denoted by
	\begin{equation}
		\label{eq:def-cal-N}
	{\cal N}={\cal N}({\cal M}, \{p_k\}_{k\in {\cal M}}, \{q_{k,j}\}_{k,j\in {\cal M}}).
	\end{equation}

The adoption level at time~$t$ is~$\frac{1}{M}\sum_{\substack{j \in {\cal M}}}X_j(t)$. Our goal is to 
obtain lower and upper bounds for
the {\bf expected adoption level} (fraction of adopters)
	\begin{equation*}
		\label{eq:number_to_fraction}
		f(t;{\cal N}):=\frac{1}{M}\mathbb{E} \Big[\sum_{\substack{j \in {\cal M}}}X_j(t)\Big].
	\end{equation*}
To do that, we will compute lower and upper bounds for  the adoption probabilities of nodes
\begin{equation*}
f_j(t;{\cal N}) := \mathbb{P}(X_j(t)=1)=\mathbb{E}\left[X_j(t)\right], \qquad j\in {\cal M},
\end{equation*}
and then use
\begin{equation}
	\label{eq:f=sum-f_j}
f= \begin{cases}
	\frac{1}{M}\sum_{j=1}^M f_j, & \quad M<\infty, \\
	\lim\limits_{M\rightarrow \infty}\frac{1}{M}\sum_{j=1}^M f_j, & \quad M = \infty.
\end{cases}
\end{equation}
The dependence of the adoption level and of the adoption probabilities of nodes on the external and internal influence rates is monotonic:
\begin{theorem}[\cite{Bass-boundary-18}]
	\label{thm:dominance-principle-f_j}
	Consider the  Bass model~\eqref{eq:dbm} on network~$\cal N$, see~\eqref{eq:def-cal-N}. Let $t>0$. Then~$f(t;{\cal N})$ is monotonically increasing, and~$\{f_m(t;{\cal N})\}$ are monotonically non-decreasing, with respect to each~$p_j$ and each~$q_{k,j}$.
\end{theorem}	

\subsection{Homogeneous complete networks}

Let $f_{\rm complete}(t;p,q,M)$ denote the expected adoption level  in the  Bass model~\eqref{eq:dbm} on the homogeneous complete network ${\cal N}_{\rm complete}(p,q,M)$, defined as 
\begin{equation}
	\label{eq:comp_net}
	p_j \equiv  p, \quad q_{k,j}=\begin{cases}
		\frac{q}{M-1}, \qquad &k\neq j,\\
		0, \qquad & k =j,
	\end{cases}, \qquad j,k\in {\cal M}.
\end{equation}
As $M$ increases, each node is influenced by more nodes, 
but the weight of each node decreases, so that the maximal rate of internal influences $q_j \equiv q$  remains unchanged, see~\eqref{eq:q_on_node}. 
Nevertheless, the expected adoption level increases with~$M$: 
\begin{lemma}[\cite{Bass-monotone-convergence-23}]
	\label{lem:f_complete-monotone-in-M}
	Let $t,p,q>0$. Then $f_{\rm complete}(t;p,q,M)$ is monotonically increasing 
	in~$M$.
\end{lemma}
As~$ M \to \infty$, the  Bass model~\eqref{eq:dbm} on complete networks
approaches the original compartmental Bass model:
\begin{theorem}[\cite{Niu-02}]
	\label{thrm:Niu}
	$
	\lim\limits_{M\rightarrow \infty} f_{\rm complete}(t;p,q,M) = f_{\rm Bass}(t;p,q),
	$
	where~$f_{\rm Bass}$ is given by~\eqref{eq:Bass_sol}.
\end{theorem}
     From Lemma~\ref{lem:f_complete-monotone-in-M} and Theorem~\ref{thrm:Niu} we have
\begin{corollary}
	\label{cor:f_complete<f_Bass}
	Let $t,p,q>0$. Then 
	\begin{equation*}
		\label{eq:f_complete<f_Bass}
		f_{\rm complete}(t;p,q,M)< f_{\rm Bass}(t;p,q), \qquad  M = 1,2, \dots
	\end{equation*} 
\end{corollary}


%
%
\section{Main results}
\label{sec:main_results}
In this section we present the main results of this paper. The proofs are given in Section~\ref{sec:proofs}.
For clarity, we formulate the results for networks  that are homogeneous in~$\{p_j\}$ and~$\{q_j\}$, i.e.,
\begin{equation}
	\label{eq:p_j,q_j-homogeneous}
	p_j\equiv p, \quad q_j \equiv q,\qquad j\in{\cal M}.
\end{equation}
This requirement can be satisfied by any graph structure that satisfies~\eqref{eq:q_j>0}, 
and not just by the complete network~\eqref{eq:comp_net}. For example, for any given network ${\cal N}({\cal M},\{p_j\}, \{q_{k,j}\})$, define network~$\widetilde{{\cal N}}({\cal M},\{\widetilde{p_j}\}, \{\widetilde{q}_{k,j}\})$ such that~$\widetilde{p_j}:=p$ and $\widetilde{q}_{k,j}:=q_{k,j}\frac{q}{q_j}$. Then~$\widetilde{{\cal N}}$ satisfies~\eqref{eq:p_j,q_j-homogeneous}, and it has the same nodes/edges structure as~$\cal N$.

The extension of the results to networks which do not satisfy~\eqref{eq:p_j,q_j-homogeneous} is discussed in Section~\ref{sec:inhomogeneous-pj-qj}. We also note that, quite often,
the difference in~$f$ between a network which is heterogeneous in~$\{p_j\}$ and~$\{q_j\}$
and the corresponding network which is homogeneous in~$\{p_j\}$ and~$\{q_j\}$
is quite small, even when the level of heterogeneity is not~\citep{OR-10,PNAS-12}.

\subsection{Non-tight universal bounds}

The following  universal lower and upper bounds are immediate:
\begin{lemma}
	  \label{lem:trivial-bounds}
		Consider the Bass model~\eqref{eq:dbm}
	on a network~${\cal N}$
	which is homogeneous in~$\{p_j\}$ and~$\{q_j\}$, see~\eqref{eq:p_j,q_j-homogeneous}. 
	Then 
	\begin{subequations}
\begin{equation}
   \label{eq:trivial-bounds-f_m}
	1-e^{-p t} \le f_m(t) \le 1-e^{-(p+q)t}, 
	 \quad t \ge 0, \quad m \in {\cal M},
\end{equation}
and so
\begin{equation}
	\label{eq:trivial-bounds-f}
1-e^{-p t} \le f(t) \le 1-e^{-(p+q)t}, 
    \quad  t \ge 0.
\end{equation}
	\end{subequations}		
\end{lemma}
\begin{pf}
Since $X_k(t) \in \{0,1 \}$ for any~$k \in \cal M$,  the adoption rate of node~$m$ is bounded by, see~\eqref{eq:lambda_j} and~\eqref{eq:q_on_node},
$$
p=p_m \le \lambda_m(t) \le p_m+\sum_{k\in {\cal M}}q_{k,m} = p+q, \quad m\in {\cal M}, \quad t \ge 0.
$$
Hence, \eqref{eq:trivial-bounds-f_m}~follows, and so 
\eqref{eq:trivial-bounds-f}~follows by~\eqref{eq:f=sum-f_j}.
%
\end{pf}

Thus, the lower and upper bounds~\eqref{eq:trivial-bounds-f_m}  for~$f_m(t)$  
correspond to 
the extreme cases when none of the other individuals adopted by time $t$,
and when all the other individuals adopted at~$t = 0+$, respectively. 
Therefore, these bounds are not expected to be tight, as indeed we will show below.

\subsection{Tight upper bound}
\label{sec:upper-bound}

If one adds edges to a network, this increases the adoption level~$f$ (Theorem~\ref{thm:dominance-principle-f_j}). 
The following two observations suggest a stronger result, namely, that
even if as we add edges, we lower the weights of the edges so as to keep $q_j \equiv q$  unchanged, the adoption level increases:
 \begin{enumerate}
 	\item The adoption level $f_{\rm complete}(t;p,q,M)$ in homogeneous complete networks is monotonically increasing in~$M$ (Lemma~\ref{lem:f_complete-monotone-in-M}).
 	
 	\item  The adoption level~$f_D(t;p,q)$ in infinite $D$-dimensional Cartesian networks, where each node is connected to its 
 	$2D$ nearest neighbors, and the weights of these edges  is~$\frac{q}{2D}$,
 	is monotonically increasing in~$D$ (this was shown numerically and asymptotically in~\cite{OR-10}).
 	  
 \end{enumerate}
 Thus,    
{\em numerous weak edges lead to a faster diffusion than a few strong ones}.
 Therefore, we can expect  that among all networks with $M$~nodes that satisfy~\eqref{eq:p_j,q_j-homogeneous}, the fastest diffusion would be on the complete network ${\cal N}_{\rm complete}(p,q,M)$, see~\eqref{eq:comp_net}, as formulated in Conjecture~\ref{conj-ub-M} below. If that is indeed the case,  then by Corollary~\ref{cor:f_complete<f_Bass}, the adoption levels 
on all networks should be bounded from above by~$f_{\rm Bass}$.
Indeed, we can rigorously prove 
\begin{theorem}
	\label{thm:bass_largest}
	Consider the Bass model~\eqref{eq:dbm}
	on a network~${\cal N}$
	which is homogeneous in~$\{p_j\}$ and~$\{q_j\}$, see~\eqref{eq:p_j,q_j-homogeneous}. 
Then
	\begin{equation}
		\label{eq:f_m<=f_Bass}
		f_{m}(t;{\cal N}) \leq f_{\rm Bass}(t;p,q), \qquad t \ge 0, \qquad m \in {\cal M},
	\end{equation}
	where~$f_{\rm Bass}$ is given by~\eqref{eq:Bass_sol}, and so
	\begin{equation}
		\label{eq:f<=f_Bass}
		f(t;{\cal N})\leq f_{\rm Bass}(t;p,q), \qquad t \ge 0.
	\end{equation}
\end{theorem}

 In Lemma~\ref{lem:trivial-bounds} we derived the upper bound
$
f_j(t), f(t) \le 1- e^{-(p+q)t}$.
The upper bound
of Theorem~\ref{thm:bass_largest} is better (i.e., lower), since by~\eqref{eq:Bass_sol},
$$
f_{\rm Bass}(t;p,q)=\frac{1-e^{-(p+q)t}}{1+\frac{q}{p}e^{-(p+q)t}}< 1- e^{-(p+q)t}.
$$ 
We can further show that $f_{\rm Bass}$ is the {\em tight} universal upper bound:
\begin{lemma}
	\label{lem:f_Bass-is-tight}
	The universal upper bound in Theorem~\ref{thm:bass_largest} is tight, in the sense that
	\begin{eqnarray*}
	&\sup_{\{ \, {\cal N} \, | \, \text{\eqref{eq:p_j,q_j-homogeneous} holds}\}} f(t;{\cal N})
	 =
	  \sup_{\{ \, {\cal N} \, | \, \text{\eqref{eq:p_j,q_j-homogeneous} holds}\}, \, m \in \cal M} f_m(t;{\cal N}) = f_{\rm Bass}(t;p,q).
	\end{eqnarray*}	
\end{lemma}
  While the upper bound~$f_{\rm Bass}$ is attained for an infinite homogeneous complete network (Theorem~\ref{thrm:Niu}), it is strict for nodes that have a finite indegree, hence for networks with a positive fraction of nodes with finite indegree:
 \begin{theorem}
	\label{thm:f<f_Bass}
Assume the conditions of Theorem~\ref{thm:bass_largest}. 
\begin{enumerate}
\item   If node~$m$ has a finite indegree,
	then 
	\begin{equation}
		\label{eq:f_m<f_Bass}
		f_m(t;{\cal N}) < f_{\rm Bass}(t;p,q), \qquad t>0.
	\end{equation}

\item If there is a positive fraction of nodes  in the network with a finite indegree, then
\begin{equation}
	 \label{eq:f<f_Bass=thm}
	f(t;{\cal N}) < f_{\rm Bass}(t;p,q), \qquad t > 0.
\end{equation}
\end{enumerate}
\end{theorem}

Therefore, the upper bound~$f_{\rm Bass}$ is strict for any network which is not  infinite and complete (up to a vanishing fraction of nodes). In particular, assume that the network type is one of the following: 
\begin{itemize}
	\item A finite network.
	\item An infinite (homogeneous or heterogeneous) D-dimensional Cartesian network.
	\item An infinite  scale-free network~\citep{barabasi1999emergence}.
	\item An infinite  small-worlds network~\citep{watts1998collective}.
	\item The infinite sparse random networks $\lim_{M \to \infty} G\left(M,\frac{\lambda}{M}\right)$ ~\citep{erdHos1960evolution}.
\end{itemize}
	Since all these finite and infinite networks have a positive fraction of finite-indegree nodes, Theorem~\ref{thm:f<f_Bass} implies that $f<f_{\rm Bass}$ for all these networks types.

\subsection{Tight lower bound}

Let ${\cal N}_{M=2}^{\rm hom}(p,q)$ denote the
 homogeneous network with two nodes, where
\begin{equation}
\label{eq:N_M=2}
{\cal M} = \{1, 2\}, \quad  p_1=p_2 = p, \quad q_{1,2} = q_{2,1} = q, \quad q_{1,1} = q_{2,2} = 0.
\end{equation}
The expected adoption level on ${\cal N}_{M=2}^{\rm hom}$ can be explicitly calculated  (see, e.g.,~\cite{Bass-boundary-18}),  giving
\begin{equation}
\label{eq:f_all-M=2}
f_{M=2}^{\rm hom}(t;p,q)  = 1-e^{-pt}\frac{qe^{-pt}-pe^{-qt}}{q-p}, \qquad p \neq q.
\end{equation}
Note that there is only one homogeneous network with two nodes. Thus,
$f_{M=2}^{\rm hom}(t;p,q)  = f_{\rm complete}^{\rm hom}(t;p,q,M=2)  = f_{\rm circle}^{\rm hom}(t;p,q,M=2)$.

The informal arguments at the beginning of Section~\ref{sec:upper-bound} 
  suggest that {\em few strong edges lead to a slower diffusion than numerous   weak ones}. Hence,  
it is intuitive to expect that for given~$p$ and~$q$, the adoption level is lowest when the influence~$q$ on any node in the network is exerted by a single node.  
   	This requirement is satisfied when the network is a one-sided circle, or a collection of disjoint one-sided circles. Among all circles, the lowest adoption is on a two-node circle (Lemma~\ref{lem:f_complete-monotone-in-M}). Intuitively, this is because 
   on a two-node circle each node can only be influenced by one node, whereas on longer circles each node can also be indirectly influenced by additional nodes. 
%
%
%
%
%
%
Indeed, we now prove that $f_{M=2}^{\rm hom}$ is a universal lower bound for~$\{f_m\}$, hence for~$f$:
\begin{theorem}
\label{thm:f_m>=f_M=2}
Assume the conditions of Theorem~\ref{thm:bass_largest}. 
Then
\begin{subequations}
	\label{eqs:f,f_m>=f_M=2}
	\begin{equation}
		\label{eq:f_m>=f_M=2}
		f_m(t;{\cal N})  \ge  f_{M=2}^{\rm hom}(t;p,q), \qquad t \ge 0, \qquad m \in {\cal M},
	\end{equation}
	and so
	\begin{equation}
		\label{eq:f>=ff_M=2}
		f(t;{\cal N})  \ge  f_{M=2}^{\rm hom}(t;p,q), \qquad t \ge 0.
	\end{equation}
\end{subequations}
\end{theorem}

In Lemma~\ref{lem:trivial-bounds}, we derived the lower bound
$
f_j(t), f(t)  \ge  1-e^{-pt}$. The lower bound
in Theorem~\ref{thm:f_m>=f_M=2} is better (i.e., larger), since by Theorem~\ref{thm:dominance-principle-f_j},
$$
f_{M=2}^{\rm hom}(t;p,q)>f_{M=2}^{\rm hom}(t;p,q=0) = 1-e^{-pt}.
$$ 

Moreover, $f_{M=2}^{\rm hom}$ is the {\em tight} universal lower bound:
\begin{lemma}
	\label{lem:f_M=2-is-tight}
	Let $M \in\{ 2,4, \dots \}$. 
	Then
	\begin{eqnarray*}
&\inf_{\{ \, {\cal N} \, | \, \text{\eqref{eq:p_j,q_j-homogeneous} holds}\}} f(t;{\cal N}) = 
  \inf_{\{ \, {\cal N} \, | \, \text{\eqref{eq:p_j,q_j-homogeneous} holds}\}, \, m \in \cal M} f_m(t;{\cal N}) = f_{M=2}^{\rm hom}(t;p,q).	
	\end{eqnarray*}
\end{lemma}

 The lower bound~$f,f_j\ge  f_{M=2}^{\rm hom}$ is attained 
  if the network is a collection of disjoint pairs of nodes, 
  each of which is of type~${\cal N}_{M=2}^{\rm hom}$. For all other networks,  however,  it is strict:
\begin{theorem}
	\label{thm:f_j>f_M=2^hom}
	Assume the conditions of Theorem~\ref{thm:bass_largest}. 
	\begin{itemize}
		\item 
		If node~$j$ belongs to a connected component with more than two nodes, then 
		\begin{equation}
			\label{eq:f_m>f_2nodes}
			f_j(t;{\cal N}) > f_{M=2}^{\rm hom}(t;p,q), \qquad t>0.
		\end{equation}
		
		\item 	If the fraction of nodes  in~$\cal N$ that belong to a connected component with more than two nodes is positive, then
	\begin{equation}
		\label{eq:f>f_M=2^hom-thm}
		f(t;{\cal N})> f_{M=2}^{\rm hom}(t;p,q), \qquad t > 0.
	\end{equation}
			\end{itemize}
\end{theorem}

\subsection{Bounds for networks   inhomogeneous in 
	$\{p_j\}$ or~$\{q_j\}$}
  \label{sec:inhomogeneous-pj-qj}

We can extend all the upper-bound results to networks which are not homogeneous in 
$\{p_j\}$ and in~$\{q_j\}$, as follows:
\begin{corollary}
	Theorem~\ref{thm:bass_largest}, Lemma~\ref{lem:f_Bass-is-tight}, 
	and  Theorem~\ref{thm:f<f_Bass},  
	remain valid if we replace  condition~\eqref{eq:p_j,q_j-homogeneous} with 
	\begin{equation}
		\label{eq:p_j,q_j<p,q}
		p_j \leq  p, \quad q_j \leq q, \qquad j\in {\cal M}.
	\end{equation} 
\end{corollary} 
\begin{pf}
This follows from Theorem~\ref{thm:dominance-principle-f_j}.
\end{pf}

Similarly, we can extend all the lower-bound results to networks which are not homogeneous in 
$\{p_j\}$ and in $\{q_j\}$: 
\begin{corollary}
	Theorem~\ref{thm:f_m>=f_M=2}, Lemma~\ref{lem:f_M=2-is-tight}, and  Theorem~\ref{thm:f_j>f_M=2^hom}, 
	remain valid if we replace condition~\eqref{eq:p_j,q_j-homogeneous} with 
	\begin{equation}
		\label{eq:p_j,q_j>=p,q}
		p_j \geq  p, \quad q_j \geq q, \qquad j\in {\cal M}.
	\end{equation}
\end{corollary}

\section{Gap between lower and upper bounds}
\label{sec:gap}

Consider any network~${\cal N}$ which is homogeneous in~$\{p_j\}$ and~$\{q_j\}$, see~\eqref{eq:p_j,q_j-homogeneous}.
By Theorems~\ref{thm:bass_largest} and~\ref{thm:f_m>=f_M=2}, the expected adoption level and the adoption probability of nodes are bounded by
  $$
f^{\rm hom}_{M=2}(t;p,q)\leq f(t;{\cal N}), f_m(t;{\cal N}) \leq f_{\rm Bass}(t;p,q),
\qquad t\ge 0.
$$  
Therefore, it is natural to consider the {\em size of the gap} between the explicit lower and upper bounds~$f_{M=2}^{\rm hom}$ and~$f_{\rm Bass}$,
which expresses the dependence of the diffusion on the network structure.

\begin{figure*}[ht!]
\begin{center}
\scalebox{1}{\includegraphics{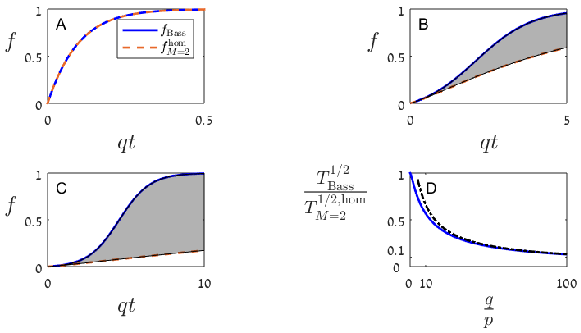}}
\end{center}
\caption{The expected adoption level~$f(t)$ of any network~${\cal N}$ lies in the shaded region between the lower bound~$f^{\rm hom}_{M=2}$ (orange dash) and upper bound~$f_{\rm Bass}$ (blue solid). 
 (A)~$\frac{q}{p}=0.1$. The two bounds are nearly indistinguishable. (B)~$\frac{q}{p}=10$. (C)~$\frac{q}{p}=100$. (D)~$\frac{T_{\rm Bass}^{1/2}}{T^{1/2, {\rm hom}}_{M=2}}$ as a function of~$\frac{q}{p}$ (blue solid), and its asymptotic approximation~\eqref{eq:approx} (black dashes). }
\label{fig:f_bass_2nodes}
\end{figure*}

The explicit bounds can be written in a dimensionless form as
  $$
f_{M=2}^{\rm hom}(t;p,q)=f_{M=2}^{\rm hom}\left(\widetilde{t};\widetilde{q}\right), \quad f_{\rm Bass}(t;p,q) = f_{\rm Bass}\left(\widetilde{t};\widetilde{q}\right),
$$
where~$\widetilde{t}=qt$ and~$\widetilde{q}=\frac{q}{p}$.
The nondimensional parameter~$\widetilde{q}$  expresses the ratio of internal and external influences. 
Since network effects are only due to internal influences, they increase with~$\frac{q}{p}$. Thus,
when~$q=0$, there are no network effects, and so the two bounds are identical, i.e., 
$$
f^{\rm hom}_{M=2}(t;p,q=0) = f_{\rm Bass}(t;p,q=0)=1-e^{-pt}.
$$
 When~$\frac{q}{p}\ll 1$,
 the network has a minor effect on the diffusion, and so~$f^{\rm hom}_{M=2}\approx f_{\rm Bass}$, see Figure~\ref{fig:f_bass_2nodes}A.  For products that spread predominantly through word-of-mouth, however, the regime of relevance is~$\frac{q}{p} \gg 1$, typically $10\leq \frac{q}{p}\leq 100$~\citep{Bass-69}. As can be expected, the difference between~$f^{\rm hom}_{M=2}$ and $f_{\rm Bass}$ is significant for~$\frac{q}{p}=10$ (Figure~\ref{fig:f_bass_2nodes}B), and even larger for~$\frac{q}{p}=100$ (Figure~\ref{fig:f_bass_2nodes}C). Note that for any network~${\cal N}$, $f(t;{\cal N})$ lies in the shaded region between~$f^{\rm hom}_{M=2}(t)$ and $f_{\rm Bass}(t)$. 

It is  instructive to compare the adoption levels on different networks using the ``half-life''~$T^{1/2}$ for half of the population to adopt.
In particular, we can use~$T^{1/2}$ to compare the bounds~$f_{\rm Bass}$ and~$f_{M=2}^{\rm hom}$.
 The ratio~$\frac{T_{\rm Bass}^{1/2}}{T^{1/2, {\rm hom}}_{M=2}}$ can be estimated asymptotically, yielding
 \begin{equation}
 \label{eq:approx}
 \frac{T_{\rm Bass}^{1/2}}{T^{1/2, {\rm hom}}_{M=2}}\sim \frac{2}{\log 2 }\frac{p}{q}\log \frac{q}{p}, \qquad  \frac{q}{p} \gg 1.
 \end{equation}
 Figure~\ref{fig:f_bass_2nodes}D confirms that~$\frac{T^{1/2}_{\rm Bass}}{T^{1/2, {\rm hom}}_{M=2}}$ decreases with~$\frac{q}{p}$, and approaches the asymptotic limit~\eqref{eq:approx} as~$\frac{q}{p}\rightarrow \infty$.
This limit goes to zero as~$\frac{q}{p}\rightarrow \infty$, showing that the network structure has a large effect on the diffusion when~$\frac{q}{p}\gg 1$, i.e., for products that  diffuse primarily by internal influences.


\section{Open problems}
  \label{sec:open}
  
This manuscript settles the conjecture from~\cite{OR-10},
 but  leads to some new questions, which are currently open. 
 Indeed, the upper and lower bounds in Theorems~\ref{thm:bass_largest}
 and~\ref{thm:f_m>=f_M=2} are tight for networks with any number of nodes. 
Can these bounds be improved if we restrict ourselves to networks with a fixed number of nodes? 

Thus, let $${\cal G}(p,q,M):=\{{\cal N}  \mid \text{${\cal N}$ has $M$ nodes, \eqref{eq:p_j,q_j-homogeneous} holds} \}$$
 be the set of all networks with $M$~nodes that are homogeneous in~$\{p_j\}$ and~$\{q_j\}$.  In the beginning  of Section~\ref{sec:upper-bound}, we argued that the fastest diffusion in~${\cal G}(p,q,M)$ should occur on the homogeneous complete network~\eqref{eq:comp_net}. Therefore, we formulate 
 \begin{conjecture}
 	\label{conj-ub-M}
 $
\sup_{{\cal N} \in {\cal G}(p,q,M)} f(t;{\cal N})  =   f_{\rm complete}(t;p,q,M).
$ 	
\end{conjecture} 	
We note, however, that the rate of convergence of~$f_{\rm complete}$ to~$f_{\rm Bass}$  as $M \to \infty$ is~$O(\frac1M)$, see~\cite{Bass-compartmental-limits-22}. Therefore, the difference between these two upper bounds becomes negligible for large (e.g., $M = 10^6$) networks.

 Consider now the lower bound.
 Let $M$ be even, and let network~$\tilde{\cal N}$ be composed of~$\frac{M}{2}$ pairs of nodes, each of which is of type~${\cal N}_{M=2}^{\rm hom}$, see~\eqref{eq:N_M=2}. Then $f(t;\tilde{\cal N}) = f_{M=2}^{\rm hom}(p,q)$. 
Therefore,
$$
\inf_{{\cal N} \in {\cal G}(p,q,M)} f(t;{\cal N})  =  f_{M=2}^{\rm hom}(p,q), 
\qquad \text{$M$ even}.
$$
Thus,  the lower bound $f_{M=2}^{\rm hom}$ cannot be improved 
(i.e., increased) for networks with a fixed even number of nodes.
The tight lower bound for $M$~odd, however, is an open problem. 

 Another open question is the tight lower bound of~$f$ among {\em connected}
networks with $M$~nodes (even or odd) that are homogeneous in~$\{p_j\}$ and~$\{q_j\}$.
Here one may need to distinguish between connected undirected networks, 
 weakly-connected directed graphs (there is an undirected path between any pair of vertices), and strongly-connected directed graphs (there is a directed path between every pair of vertices).


\section{Proof of results}
\label{sec:proofs}

\subsection{Master equations}

Denote the nonadoption probability of node~$j$ by
\begin{equation}
\label{eq:nonadoption_prob}
[S_j](t):=1-f_j(t)= \mathbb{P}(X_j(t)=0).
\end{equation}
Then $[S_j]$ satisfies
 the {\em master equation}~\citep{fibich2022diffusion}
	\begin{equation}
		\label{eq:master_Sa}
		\frac{d}{dt}[S_j](t)= -\left(p_j+q_j\right)\left[S_{j}\right]
		+\sum_{\substack{k \in {\cal M}}}q_{k,j}[S_j,S_k](t), \qquad [S_{j}](0) =1,
	\end{equation} 
where~$q_j$ is given by~\eqref{eq:q_on_node}, and
$$
[S_j,S_k](t):=\mathbb{P}(X_j(t)=X_k(t)=0).
$$
In general, to close these equations, one adds the master equations for all pairs~$\{[S_j, S_k]\}$, all triplets~$\{[S_j, S_k, S_m]\}$, etc., see~\cite{fibich2022diffusion}. For the purpose of obtaining the lower and upper bounds, however, we will only need the following result:
\begin{lemma}
Consider the  Bass model~\eqref{eq:dbm}. Then for any~$i,j \in {\cal M}$,
\begin{equation}
\label{eq:set_e-2pt}
[S_i](t)\,[S_j](t) \leq [S_i,S_j](t) \leq e^{-2pt}, \qquad 0\leq t<\infty.
\end{equation}
\end{lemma}
\begin{pf}
The left inequality is proved in~\cite{Funnel2022}. For the right inequality, we note that the joint nonadoption probability of a pair~$\{i,j\}$ of isolated nodes ($q_{j} = q_{i} = 0$) is given by
		\begin{equation*}
			[S_i,S_j](t;p,q_{j} = q_{i} = 0) = e^{-2pt}, 
		\end{equation*}
see, e.g.,~\cite{Bass-boundary-18}. Hence, the right inequality follows from Theorem~\ref{thm:dominance-principle-f_j}.
\end{pf}
\subsection{Differential and integral Bass inequalities}
Let us recall the following result:
\begin{lemma}[\cite{OR-10}]
	\label{lem:Bass-inequality}
	Let~$p,q>0$, and let $f(t)$ satisfy the {\bf differential Bass inequality}
	$$
	\frac{df}{dt} < (1-f)(p+qf), \quad t>0, \qquad f(0)=0. 
	$$
	Then $f(t) < f_{\rm Bass}(t;p,q)$ for $0<t<\infty$.
\end{lemma}

Let $[S_{\rm Bass}]:=1-f_{\rm Bass}$ denote the nonadoption level in the compartmental
Bass model. Then by~\eqref{eq:homogeneous_Bass_frac},
\begin{equation}
	\label{eq:homogeneous_Bass_frac-[S]}
	\frac{d}{dt}[S_{\rm Bass}](t)
	 =-(p+q)[S_{\rm Bass}]+ q[S_{\rm Bass}]^2  , \quad [S_{\rm Bass}](0)=1.
\end{equation}
If we replace the equality sign in~\eqref{eq:homogeneous_Bass_frac-[S]} by an inequality, the solution of this inequality is bounded from below by~$[S_{\rm Bass}]$:
 \begin{lemma}
 \label{cor:diffBass}
Let~$p,q>0$, and let $[S](t)$ satisfies the {\bf differential Bass inequality}
\begin{equation*}
	\label{eq:homogeneous_Bass_frac-[S]-inequality}
	\frac{d}{dt}[S](t)
	>-(p+q)[S]+ q[S]^2, \quad t>0, \qquad [S](0)=1.
\end{equation*}
 Then $[S](t)> [S_{\rm Bass}](t)$ for $0<t<\infty$.
 \end{lemma}
\begin{pf}
This follows from Lemma~\ref{lem:Bass-inequality}  and~$[S_{\rm Bass}]=1-f_{\rm Bass}$.
\end{pf}

Multiplying~\eqref{eq:homogeneous_Bass_frac-[S]} by $e^{(p+q)t}$, integrating between zero and~$t$, and using the initial condition, gives the {\em integral form of the compartmental Bass model}
\begin{equation}
	\label{eq:homogeneous_Bass_frac-S-integral}
	[S_{\rm Bass}] (t) = e^{-(p+q)t} + q\int_0^t e^{-(p+q)(t-\tau)} [S_{\rm Bass}]^2(\tau) \, d\tau.
\end{equation}
If we replace the
equality sign in~\eqref{eq:homogeneous_Bass_frac-S-integral} by an inequality, the solution of the resulting  integral Bass inequality is bounded from below by~$[S_{\rm Bass}]$:
\begin{lemma}
	\label{lem:Bass-inequality-integral}
	Let $p,q>0$, and let $[S](t)$ be non-negative and continuous in~$[0, \infty)$. 
	\begin{enumerate}
		\item 
	If $[S]$ satisfies the {\bf  integral Bass inequality}
\begin{equation}
	\label{eq:S-Bass-inequality-integral}
	[S] (t) \ge  e^{-(p+q)t} + q\int_0^t e^{-(p+q)(t-\tau)} [S]^2(\tau) \, d\tau, \quad t>0,
\end{equation}
	then 
	$
	[S] (t) \ge 	[S_{\rm Bass}] (t;p,q)$ for $t \ge 0$.

	\item If  inequality~\eqref{eq:S-Bass-inequality-integral} is strict, then  
	$[S] (t) > 	[S_{\rm Bass}] (t;p,q)$ for $t>0$.
		\end{enumerate}
\end{lemma}
\begin{pf}
Let $u:=[S]-[S_{\rm Bass}]$. Subtracting~\eqref{eq:homogeneous_Bass_frac-S-integral} from~\eqref{eq:S-Bass-inequality-integral} gives
$$
  u(t) \ge q \int_0^t e^{-(p+q)(t-\tau)}  \left([S]^2-[S_{\rm Bass}]^2\right)(\tau) \, d\tau.
$$
Therefore,
\begin{eqnarray}
	\label{eq:u>=int-phi*u}
	u(t) \ge \int_0^t\phi(\tau) u(\tau) \, d\tau, 
	 \qquad \phi(\tau):=q e^{-(p+q)(t-\tau)}  \left([S]+[S_{\rm Bass}]\right)(\tau).
\end{eqnarray}
%
%
Since $[S]$ and $[S_{\rm Bass}]$ are continuous and non-negative, then so is~$\phi$.
Let 
\begin{equation}
	\label{eq:v-def-integral-Bass}
 v(t):=e^{-\int_0^t \phi} \int_0^t\phi(\tau) u(\tau) \, d\tau.
\end{equation}
 Then
$v(0) = 0$  and
$$
  \frac{dv}{dt} = 
  e^{-\int_0^t \phi} \phi(t) 
  \left(u(t)- \int_0^t\phi(\tau) u(\tau) \, d\tau \right) \ge 0,
$$
where the inequality follows from~\eqref{eq:u>=int-phi*u}. Therefore,  for $t \ge 0$,  $v(t) \ge 0$.
Hence, by~\eqref{eq:v-def-integral-Bass}, $\int_0^t\phi(\tau) u(\tau) \, d\tau \ge 0$ and so by~\eqref{eq:u>=int-phi*u}, $u(t) \ge 0$.

If inequality~\eqref{eq:S-Bass-inequality-integral} is strict, we replace in the above proof all
 ``$\ge$'' signs by~``$>$'' signs.
\end{pf}

\subsection{Upper bound}
 \label{sec:upper}

We begin with an auxiliary result.
\begin{lemma}
	\label{cor:lipschitz}
	Consider the  Bass model~\eqref{eq:dbm}.
	Let~\eqref{eq:p_j,q_j<p,q} hold, and let
	\begin{equation}
	\label{eq:s_inf_s_j}
		[\underbar{S}](t) :=  \inf_{j\in {\cal M}} \left\{[S_j](t) \right\}.
	\end{equation}
	Then $[\underbar{S}](t)$ is non-negative and continuous.
\end{lemma}
\begin{pf}
	The non-negativity of~$[\underbar{S}]$ follows from that of $\{[S_j]\}$. 
	Let $j \in {\cal M}$.
	Since all probabilities are bounded between~$0$ and~$1$, then using~\eqref{eq:master_Sa} and~\eqref{eq:p_j,q_j<p,q},
	\begin{equation*}
		\left|\frac{d}{dt}[S_j]\right| \leq (p+q)[S_j]+\sum_{k\in \cal M}q_{k,j}\left[S_j,S_k\right]\leq p+q+\sum_{k\in \cal M}q_{k,j} \le  \kappa, 
	\end{equation*}
	where~$\kappa:=p+2q$.
	Hence, by the mean-value theorem, for any $t,t^*>0$, $\left|[S_j](t) - [S_j]\left(t^*\right)\right|\leq \kappa\left|t-t^*\right|$, and so~$- [S_j]\left(t^*\right)\leq -[S_j](t)+\kappa\left|t-t^*\right|\leq -[\underbar{S}](t)+\kappa\left|t-t^*\right|$. Taking the supremum of the left-hand side yields $ -[\underbar{S}]\left(t^*\right)\leq -[\underbar{S}](t)+\kappa\left|t-t^*\right|$, and so~$[\underbar{S}](t) -[\underbar{S}]\left(t^*\right)\leq \kappa\left|t-t^*\right|$. Swapping~$t$ and~$t^*$ gives the inverse estimate, and so~$[\underbar{S}](t)$ is continuous.
\end{pf}

{\bf Proof of Theorem~\ref{thm:bass_largest}. }
Since~$1-f_m = [S_m]\geq [\underbar{S}]$, 
see~\eqref{eq:nonadoption_prob} and~\eqref{eq:s_inf_s_j}, 	
	it is sufficient to show that
	\begin{equation}
		\label{eq:s_bass<s_inf}
		[\underbar{S}](t) \geq [S_{\rm Bass}](t,p,q).
	\end{equation}
	
	By~\eqref{eq:master_Sa} with $q_j = q$, see~\eqref{eq:p_j,q_j-homogeneous}, 
	\begin{equation}
		\label{eq:S_j-imtegral-identity-upper-bound}
		[S_j] = e^{-(p+q)t}+\int_0^te^{-(p+q)\left(t-\tau\right)}
		\sum_{k\in \cal M} q_{k,j} \, [S_j,S_k](\tau)\, d\tau .
	\end{equation}
	Therefore,  by the lower bound in~\eqref{eq:set_e-2pt} and~\eqref{eq:s_inf_s_j},
	\begin{eqnarray*}
			&[S_j] \geq e^{-(p+q)t}+\int_0^te^{-(p+q)(t-\tau)}\sum\limits_{k\in \cal M}q_{k,j} \, [S_j](\tau) \, [S_k](\tau) \, d\tau \\
			&\geq e^{-(p+q)t}+q \int_0^te^{-(p+q)(t-\tau)} [\underbar{S}]^2(\tau) d\tau.
	\end{eqnarray*}
	Taking the infimum over all~$j$ gives
	\begin{equation*}
		[\underbar{S}] \geq e^{-(p+q)t}+q \int_0^t e^{-(p+q)(t-\tau)}[\underbar{S}]^2(\tau)\, d\tau.
	\end{equation*}
	Therefore, since $[\underbar{S}]$ is non-negative and continuous (Lemma~\ref{cor:lipschitz}),
	we can use the integral Bass inequality (Lemma~\ref{lem:Bass-inequality-integral})
	to get inequality~\eqref{eq:s_bass<s_inf}, from which~\eqref{eq:f_m<=f_Bass} follows. Therefore, by~\eqref{eq:f=sum-f_j}, 
\eqref{eq:f<=f_Bass}~follows.~\qed

{\bf Proof of Lemma~\ref{lem:f_Bass-is-tight}. }
The result for~$f$ follows from  Theorem~\ref{thrm:Niu}.  Since the complete network~\eqref{eq:comp_net} is homogeneous, 
	$f_m \equiv f$ for all $m \in \cal M$. Hence, the result holds for any~$f_m$ as well. \qed

{\bf Proof of Theorem~\ref{thm:f<f_Bass}. }
	Let 
$$
A_{d}({\cal N}):= \{m \in {\cal M} \mid \text{indegree}\left(m\right)=d \}
$$
denote  the set of all nodes with
indegree~$d$ in network~${\cal N}$. Then it is sufficient to prove that for all networks 
that satisfy~\eqref{eq:p_j,q_j-homogeneous}  and for all $d \in \mathbb{N}$, 
\begin{equation}
	\label{eq:S_m>S_Bass-induction}
	[S_m](t;{\cal N}) > [S_{\rm Bass}](t;p,q), \quad t>0, \qquad m \in A_d({\cal N}).
\end{equation}

We prove~\eqref{eq:S_m>S_Bass-induction} by induction on~$d$. 
When $d=0$, node $m \in A_0$ is not influenced by any other node, and so 
\begin{equation}
	\label{eq:e^-pt>S_Bass}
	[S_m](t;{\cal N})= e^{-pt}   =[S_{\rm Bass}](t;p,q=0)
	> [S_{\rm Bass}](t;p,q),
\end{equation}
where the inequality follows from Theorem~\ref{thm:dominance-principle-f_j}.

For the induction stage, we assume that~\eqref{eq:S_m>S_Bass-induction}
holds  for all networks 
that satisfy~\eqref{eq:p_j,q_j-homogeneous}  and for all~$m \in A_{d-1}$, and prove that it holds 
for all networks that satisfy~\eqref{eq:p_j,q_j-homogeneous}  and for all~$m \in A_{d}$, as follows. 
Let $m \in A_d$, where $d \ge 1$, and denote by 
$\{k_1,\ldots, k_{d}\}$ the $d$~nodes that can influence~$m$. 
The master equation for~$[S_m]$ is,
see~\eqref{eq:p_j,q_j-homogeneous} and~\eqref{eq:master_Sa}, 
\begin{equation}
	\label{eq:s^j_evolution_J+1}
	\frac{d}{dt}[S_m]= -(p+q)[S_m]
	+\sum_{i=1}^{d} q_{k_i,m}[S_m,S_{k_i}], \qquad [S_m](0)=1 .
\end{equation}
By the indifference principle, we can compute each of the $d$~probabilities $\{[S_m,S_{k_i}]\}_{i=1}^{d}$  on a modified network~$\widetilde{\cal N}_i$, in which we remove the edge~$k_i \rightarrow {m}$. Thus, $	[S_m,S_{k_i}] = \widetilde{[S_m,S_{k_i}]}$,  where the tilde sign refers to probabilities in~$\widetilde{\cal N}_i$. In this modified network, node~$m$ has indegree~$d-1$, and so by the induction assumption\footnote{In the modified network~$\widetilde{\cal N}_i$ we reduced~$q_m$ by~$q_{k_i,m}>0$. Therefore, $\widetilde{q}_m<q$,
	and so we cannot apply the induction assumption directly for~$\widetilde{\cal N}_i$.  
	By Theorem~\ref{thm:dominance-principle-f_j}, however,
	aince the induction assumption holds when~$\widetilde{q}_m=q$, 
	see~\eqref{eq:p_j,q_j-homogeneous}, it also holds when~$\widetilde{q}_m<q$.}
$$
\widetilde{[{S}_m]} > [S_{\rm Bass}]. 
$$
In addition, by Theorem~\ref{thm:bass_largest}, 
$$
\widetilde{[{S}_{k_i}]}  \ge  [S_{\rm Bass}], 
$$
Combining the above and~\eqref{eq:set_e-2pt},
we have that
$$
[S_m,S_{k_i}] = \widetilde{[S_m,S_{k_i}]} \ge \widetilde{[{S}_m]}\widetilde{[{S}_{k_i}]} >[S_{\rm Bass}]^2.
$$
Therefore, 
\begin{equation}
	\label{eq:s^js^i_indifference_J+1}
	\sum_{i=1}^{d} q_{k_i,m} \, [S_m,S_{k_i}] > \sum_{i=1}^{d} q_{k_i,m} [S_{\rm Bass}]^2 = q [S_{\rm Bass}]^2.
\end{equation}

By~\eqref{eq:s^j_evolution_J+1} and~\eqref{eq:s^js^i_indifference_J+1},
\begin{equation*}
	\frac{d}{dt}[S_m]+(p+q)[S_m]>q[S_{\rm Bass}]^2, \qquad [S_m](0)=1 .
\end{equation*}
This is the differential Bass inequality (Lemma~\ref{lem:Bass-inequality}), written in terms
of~$[S]$, see~\eqref{eq:homogeneous_Bass_frac-[S]-inequality}.  Hence, 
$[S_m] > [S_{\rm Bass}]$, as needed.
\qed

\subsection{Lower bound}

{\bf Proof of Theorem~\ref{thm:f_m>=f_M=2}. }
	To prove the lower bound~\eqref{eq:f_m>=f_M=2} for~$f_m$,
	it is sufficient to show that
	\begin{eqnarray*}
		[S_m](t) \leq [S_{M=2}^{\rm hom}](t;p,q) :=1-f_{M=2}^{\rm hom}(t;p,q) 
		=  e^{-pt}\frac{qe^{-pt}-pe^{-qt}}{q-p},
	\end{eqnarray*}
where $[S_m] = 1-f_m$. 
	By the upper bound in~\eqref{eq:set_e-2pt} and~\eqref{eq:S_j-imtegral-identity-upper-bound}, we have that
		\begin{equation*}
			\begin{aligned}
			 & [S_m]  \leq e^{-(p+q)t}+\int_0^te^{-(p+q)\left(t-\tau\right)}\sum_{k \in \cal M}q_{k,m}e^{-2p\tau}\, d\tau
			= e^{-(p+q)t}+q\int_0^te^{-(p+q)\left(t-\tau\right)}e^{-2p\tau}\, d\tau 
			\\
			&~~~= \left(1-\frac{q}{q-p}\right)e^{-(p+q)t}+\frac{q}{q-p}e^{-2pt} = \left[S_{M=2}^{\rm hom}\right](t;p,q).
		\end{aligned}
	\end{equation*}
	Therefore, we proved~\eqref{eq:f_m>=f_M=2}, which implies~\eqref{eq:f>=ff_M=2}. \qed

{\bf Proof of Theoren~\ref{thm:f_j>f_M=2^hom}}
	The only inequality in the proof of Theorem~\ref{thm:f_m>=f_M=2}
	arises from using the  upper bound in~\eqref{eq:set_e-2pt}.
	Therefore,  the lower bound~\eqref{eq:f_m>=f_M=2} for~$[S_m]$ becomes an equality if and only if  $	[S_m,S_k] = e^{-2 pt}$ for all $k \in {\cal M} \setminus m$ for which $q_{k,m}>0$. A minor modification of
	Theorem~\ref{thm:dominance-principle-f_j} shows that  
	$$
	[S_j,S_k] = e^{-2 pt} \quad \iff \quad \text{$j$ and $k$ are not influenced by any other node}.
	$$
	Therefore, \eqref{eq:f_m>f_2nodes}~follows.
%
	Since $f = \frac{1}{M}\sum_{j=1}^M f_j$ for finite networks and
	$f = \lim_{M\rightarrow \infty }\frac{1}{M}\sum_{j=1}^M f_j$ for infinite networks, 	\eqref{eq:f>f_M=2^hom-thm} also follows. \qed

{\bf Proof of Lemma~\ref{lem:f_M=2-is-tight}. }
	When~$M=2$, this bound is attained by ${\cal N}={\cal N}_{M=2}^{\rm hom}(p,q)$. Moreover, 
	this bound is also attained by any finite or infinite network which is 
	a collection of disjoint pairs of nodes, each of which is of type~${\cal N}_{M=2}^{\rm hom}(p,q)$. \qed

\subsection{Asymptotic evaluation of~$\frac{T_{\rm Bass}^{1/2}}{T^{1/2, {\rm hom}}_{M=2}}$}
By~\eqref{eq:f_all-M=2}, $T^{1/2}:=T_{M=2}^{1/2, {\rm hom}}$ is the solution of 
\begin{equation}
\label{eq:t_approx}
e^{-pT^{1/2}}\frac{qe^{-pT^{1/2}}-pe^{-qT^{1/2}}}{q-p} = \frac{1}{2}.
\end{equation}
Let~$X:=e^{-pT^{1/2}}$ and~$\lambda:=\frac{q}{p}$. Then~$e^{-qT^{1/2}}=X^{\frac{q}{p}}=X^{\lambda}$. Plugging this into~\eqref{eq:t_approx}, and noting that~$0<X<1$ and~$\lambda > 0$ gives
\begin{equation*}
\label{eq:X_first}
X^2 - \frac{1}{2} = \frac{p}{q-p}\left(-X^2+X^{\lambda}\right)=O\left(\frac{1}{\lambda}\right), \qquad \lambda\gg 1.
\end{equation*}  
Therefore, 
\begin{equation*}
\label{eq:X_1^2}
X^2 \sim \frac{1}{2}, \qquad \lambda\gg 1.
\end{equation*}
Hence, by the definition of~$X$
\begin{equation*}
\label{eq:T_X}
T^{1/2, {\rm hom}}_{M=2} = \frac{1}{2p}\log(X^{-2})\sim \frac{\log(2)}{2p}, \qquad \lambda\gg 1.
\end{equation*}
Finally, by~\citealp[Lemma 11]{OR-10}, 
$$T_{\rm Bass}^{1/2} = \frac{\log\left(2+\frac{q}{p}\right)}{p+q}\sim \frac{\log\left(\frac{q}{p}\right)}{q},  \qquad \lambda\gg 1,
$$
 and so~\eqref{eq:approx} follows. \qed

\section*{Acknowledgements}
We thank Eilon Solan for useful comments.



  \bibliographystyle{apalike} 
  \bibliography{diffusion.bib}


\end{document}